\documentclass[10pt]{article}
\begin{document}

\newcommand{\nc}{\newcommand}

\newtheorem{lemma}{Lemma}[section]
\newtheorem{theorem}[lemma]{Theorem}
\newtheorem{proposition}[lemma]{Proposition}
\newtheorem{corollary}[lemma]{Corollary}
\newtheorem{remark}[lemma]{Remark}
\newtheorem{example}[lemma]{Example}
\newtheorem{hypothesis}[lemma]{Hypothesis}
\newtheorem{notation}[lemma]{Notation}
\newtheorem{definition}[lemma]{Definition}
\newtheorem{conclusion}[lemma]{Conclusion}

\nc{\QED}{\mbox{}\hfill \raisebox{-0.2pt}{\rule{5.6pt}{6pt}\rule{0pt}{0pt}}
          \medskip\par}
\newenvironment{Proof}{\noindent
    \parindent=0pt\abovedisplayskip = 0.5\abovedisplayskip
    \belowdisplayskip=\abovedisplayskip{\bf Proof. }}{\QED}
\newenvironment{Proofof}[1]{\noindent
    \parindent=0pt\abovedisplayskip = 0.5\abovedisplayskip
    \belowdisplayskip=\abovedisplayskip{\bf Proof of #1. }}{\QED}
\newenvironment{Example}{\begin{example}
                \parindent=0pt \rm}{\QED\end{example}}
\newenvironment{Remark}{\begin{remark}
                \parindent=0pt\rm }{\QED\end{remark}}
\newenvironment{Definition}{\begin{definition}
                \parindent=0pt \rm}{\QED\end{definition}}
\newenvironment{Conclusion}{\begin{conclusion}
                \parindent=0pt \rm}{\QED\end{conclusion}}

\def\ba{\begin{array}}
\def\ea{\end{array}}

\def\be{\begin{equation}}
\def\ee{\end{equation}}
\def\vs5{\vspace{0.5cm}}
\def\lab{\label}
\def\bthm{\begin{theorem}}
\def\ethm{\end{theorem}}
\def\bP{\begin{Proof}}
\def\eP{\end{Proof}}
\def\brem{\begin{remark}}
\def\erem{\end{remark}}
\def\bex{\begin{example}}
\def\eex{\end{example}}
\def\bcor{\begin{corollary}}
\def\ecor{\end{corollary}}
\def\bd{\begin{definition}}
\def\ed{\end{definition}}
\def\bprop{\begin{proposition}}
\def\eprop{\end{proposition}}
\def\blem{\begin{lemma}}
\def\elem{\end{lemma}}
\def\beq{\begin{eqnarray}}
\def\eeq{\end{eqnarray}}
\def\beqs{\begin{eqnarray*}}
\def\eeqs{\end{eqnarray*}}

\def\a{\alpha}
\def\b{\beta}
\def\s{\sigma}
\def\Sig{\Sigma}
\def\d{\delta}
\def\l{\lambda}
\def\hs{\hat{\sigma}}
\def\hT{\hat{T}}
\def\hU{\hat{U}}
\def\eps{\epsilon}
\def\veps{\varepsilon}
\def\benum{\begin{enumerate}}
\def\eenum{\end{enumerate}}
\def\bit{\begin{itemize}}
\def\eit{\end{itemize}}
\def\la{\langle}
\def\ra{\rangle}
\def\wto{\rightharpoonup}
\def\bwto{\buildrel{w(\mu_h,\mu)}\over\longrightarrow}
\def\bYto{\buildrel{Y(\mu_h,\mu)}\over\longrightarrow}
\def\bsto{\buildrel{s(\mu_h,\mu)}\over\longrightarrow}
\def\bswto{\buildrel{sw(\mu_h,\mu)}\over\longrightarrow}

\def\Cinf{C^\infty}
\def\Linf{L^\infty}
\def\supp{{\rm supp}}
\def\sm{\setminus}
\def\id{{\rm id}}

\def\dis{\displaystyle}

\def\R{{\bf R}}
\def\N{{\rm N}}
\def\Ra{{\rm R}}
\def\M{{\cal M}}
\def\B{{\cal B}}
\def\A{{\cal A}}
\def\C{{\cal C}}
\def\D{{\bf D}}
\def\calS{{\cal S}}
\def\O{{\cal O}}
\def\P{{\cal P}}
\def\Q{{\cal Q}}
\def\u{{\cal U}}
\def\V{{\cal V}}
\def\F{{\cal F}}
\def\G{{\cal G}}
\def\H{{\cal H}}
\def\E{{\cal E}}
\def\m{{\bf m}}
\def\o{{\bf o}}
\def\LL{{\cal L}}
\def\MM{{\bf M}}
\def\DD{{\cal D}}
\def\RR{{\cal R}}
\def\NN{{\bf N}}
\def\EE{{\bf E}}
\def\FF{{\bf F}}
\def\GG{{\bf G}}
\def\T{{\cal T}}
\def\TT{{\bf T}}
\def\BB{{\bf B}}
\def\HH{{\bf H}}
\def\KP{{\bf KP}}
\def\K{{\bf K}}
\def\aa{{\bf a}}

\def\rank{{\rm Rank}}
\def\span{{\rm span}}
\def\dim{{\rm dim}}
\def\diam{{\rm diam}}
\def\diag{{\rm diag}}
\def\dist{{\rm dist}}
\def\trace{{\rm trace}}
\def\div{{\rm div}}
\def\spt{{\rm spt}}
\def\Lip{{\rm Lip}}
\def\Sp{{\bf S}}
\def\Z{{\bf Z}}

\def\Grass{{\rm Grass}}
\def\diam{{\rm diam}\,}
\def\rmB{{\rm B}}

\def\pa{\partial}
\def\um{^{-1}}
\def\pim{\pi^{-1}}
\def\pm{p^{-1}}
\def\cit{\cite{BouchitteBF2}}
\def\citd{\cite{GiaquintaM2}}
\def\ov{\overline}
\def\na{\nabla}
\def\vpi{\varpi}

\def\ti{\tilde}
\def\ess{{\rm ess}}
\def\ext{^{\rm ext}}
\def\u{{\bf u}}

\def\w{{\bf w}}
\def\e{{\bf e}}
\def\Om{\Omega}
\def\om{\omega}
\def\al{\aleph}

\def\ioe{\int_{\Omega_\veps}}
\def\ioesmbe{\int_{\Omega_\veps\setminus \C_\veps}}
\def\io{\int_{\Omega}}
\def\ise{\int_{T_\veps}}
\def\ibe{\int_{\C\veps}}
\def\iyke{\int_{Y^k_\veps}}
\def\dme{\;dm_\veps}
\def\Ome{{\Omega_\veps}}
\def\me{{m_{\veps}}}
\def\intb{{\int\!\!\!\!\!\!-}}

\title{Asymptotics of a thermal flow with highly conductive
and radiant suspensions}

\author{Fadila~Bentalha $^*$, Isabelle~Gruais $^{**}$ and
Dan~Poli\v{s}evski $^{***}$
}
\date{}
\maketitle

{\bf Abstract.} Radiant spherical suspensions have an
$\veps$-periodic distribution in a tridimensional incompressible
viscous fluid governed by the Stokes-Boussinesq system. We perform
the homogenization procedure when the radius of the solid spheres
is of order $\veps^3$ (the critical size of perforations for the
Navier-Stokes system) and when the ratio of the fluid/solid
conductivities is of order $\veps^6$, the order of the total
volume of suspensions. Adapting the methods used in the study of
small inclusions, we prove that the macroscopic behavior is
described by a Brinkman-Boussinesq type law and two coupled heat
equations, where certain capacities of the suspensions and of the
radiant sources appear.



\section{Preliminaries}

One main achievement of homogenization theory Êwas the
ability to conceptually clarify Êthe relationship between microscopic and
macroscopic properties of physical systems, at least as far as the
periodic approximation could be acceptable. The major restriction was the
technically impossible interplay between different scales: if some
quantity varies as the power $\veps^\a$ of the size $\veps$ of the mesh,
then the case where $\a<0$ leads to blow up at the limit. 
This type of problems were
introduced and solved for the first time by
\cite{CioranescuM3} and developed by
\cite{BellieudB,CasadoDiaz,BrianeT,BellieudG,Mosco}.
One major
contribution in that direction is the paper by G. Allaire \cite{Allaire}
who clearly underlies the role of critical discriminating scales beyond
which nothing can be said, Êbut rigidification of elastic systems
for instance, and that can however generate a transition state where
either 'non local' effects \cite{BellieudB,BellieudG} or 'coming from
nowhere' terms \cite{CioranescuM3} can emerge. 

In this paper, we are insterested in the former case which has been
thoroughly explored when non local effects concentrate on rod-like
one-dimensional submanifolds of the three-dimensional space: see
\cite{BellieudB} for the Laplacian, \cite{BellieudG} for the Elasticity
system. This geometry enables the formulation of the limit problem as a
rod-like boundary value problem solved by the density of a Radon measure.
Our question then was: what happens in other geometries, especially if
non local effects are to be supported by a cloud of little particles? The
physical opportunity was the example of thermal flows (see
\cite{EneP,Polisevski}) where highly heat conducting spheres are immerged
in a Stokes-Boussinesq fluid. It is straightforward that for some
critical size of the particles (eventually
$\veps^3$ when the period of the distribution is $\veps$) the resulting
mixture will display a specific behaviour strongly discriminating between
a trivial case and a classically homogenized case. Our concern was 
then to develop new skills to understand how the expected Ênon local
effects would be formulated. We found out that the Dirac structure of
the masses make the classical formulation in terms of a jump term
updated Êand that it rather generates an
additional source coupled with a capacitary term representative of a
Brinkman-Boussinesq type law. Ê

More precisely, the physics of the problem may be described as follows. 
Solid spherical suspensions are $\veps$-periodically distributed in a
tridimensional bounded domain filled with an incompressible fluid
governed by the Stokes-Boussinesq system. We study the homogenization of
the convective movement which is generated by highly heterogeneous
radiant sources, when the radius of the suspensions is of 
$\veps^3$-order, that is the border case for the Navier-Stokes system
(see \cite{Allaire}). Assuming that the conductivity and the radiant
source of the fluid have $\veps^0$-order, we found that the only regular
case in which we have macroscopic effects from both the conductivity and
the radiation of the suspensions is when they are of 
$\veps^6$-order. Therefore, we have treated here strictly this case.
Nevertheless, the present procedure can be easily adapted to the other
cases.

Let $\Om\subset\R^3$ be a bounded open set and let
$$
  Y:= \left(- \frac{1}2,+\frac{1}2\right)^3.
$$
$$
  Y^k_\veps := \veps k + \veps Y,\quad k\in\Z^3.
$$
$$
  \Z_\veps := \{k\in\Z^3,\quad Y^k_\veps\subset \Om\}
$$

The reunion of the suspensions is defined by
$$
T_\veps := \cup_{k\in\Z_\veps}B(\veps k,r_\veps),
$$
where $0 < r_\veps <<\veps$ and
$B(\veps k,r_\veps)$ is the ball of radius $r_\veps$ centered at
$\veps k$, $k\in\Z_\veps$.

The fluid domain is given by
$$
  \Ome = \Om\sm T_\veps.
$$

Let $\e^{(3)}$ the last vector of the
canonical basis of
$\R^3$, $n$ the normal on $\pa\Ome$ in the outward direction and
$[\cdot]_\veps$ the jump across the interface $\pa T_\veps$.

For $a>0$ (the so-called Rayleigh number), $b>0$ (
$
  b\left(\frac{\veps}{r_\veps}\right)^3
$
denoting the ratio of the solid/fluid conductivities),  $f\in C_c(\Om)$,
$g\in C_c(\Om)$, where
$$
  C_c(\Om) := \{g\in C(\Om);\quad\supp g\quad\mbox{is compact}\quad\},
$$
we  consider the problem corresponding to
the non-dimensional Stokes-Boussinesq system governing the thermal flow
of an
$\veps$-periodic distribution  suspension of solid spheres:

To find $(u^\veps,p^\veps)$, $\theta^\veps, \zeta^\veps$
solution  of
\beq
\lab{pbe1}
\div u^\veps & = & 0 ,\quad\mbox{in}\quad \Om_\veps,
\\
- \Delta u^\veps + \na p^\veps & = & a
\theta^\veps\e^{(3)},\quad\mbox{in}\quad
\Om_\veps,
\\
- \Delta\theta^\veps + u^\veps\na\theta^\veps & = &
f,\quad\mbox{in}\quad\Om_\veps,
\\
- \Delta\zeta^\veps & = &  g,\quad\mbox{in}\quad T_\veps,
\\
\lab{zte}
\zeta^\veps &=& \theta^\veps,\quad\mbox{on}\quad\pa T_\veps
\\
\frac{\pa\theta^\veps}{\pa n}
& = & b\left(\frac{\veps}{r_\veps}\right)^3
\frac{\pa\zeta^\veps}{\pa n},\quad\mbox{on}\quad\pa T_\veps
\\
u^\veps & = & 0,\quad\mbox{on}\quad\pa\Omega_\veps,
\\
\lab{pbe8}
\theta^\veps & = & 0,\quad\mbox{on}\quad\pa\Omega.
\eeq

Set
$$
  V_\veps:= \{v\in H^1_0(\Om_\veps;\R^3),\quad \div \,v = 0\}.
$$

Thanks to (\ref{zte}),  we extend $\theta^\veps$ on $T_\veps$ by setting
$$
  \theta^\veps = \zeta^\veps\quad\mbox{on}\quad T_\veps.
$$

Then, the variational formulation reads:
\be
\lab{vare1}
\ba{lll}
\dis\forall (v,q)\in V_\eps\times L^2(\Om_\veps),&
\dis\ioe\na u^\veps\cdot\na v\;dx & =
\dis a\ioe\theta^\veps v_3\;dx
\\
\quad &
\dis\ioe q \,\div u^\veps \;dx &=  0
\ea
\ee
\be
\lab{vare4}
\ba{lcc}
\dis\forall \varphi\in H^1_0(\Om_\veps),&
\dis\ioe\!\!\na\theta^\veps\na\varphi\;dx
& +\quad
\dis b\left(\frac{\veps}{r_\veps}\right)^3
\dis\ise\na\theta^\veps\na\varphi\;dx
\\
\quad & +
\dis\ioe u^\veps\varphi\na\theta^\veps\;dx
\!\!\!
&= \!\!\!\dis\ioe\!\!\!\!f\varphi dx +
\dis b\left(\frac{\veps}{r_\veps}\right)^3\!\!\!
\dis\ise  g\varphi dx.
\ea
\ee


We define $\F_\veps\in H^{-1}(\Om)$ by
\be
\lab{fe}
 \quad\forall\varphi\in H^1_0(\Om),\quad\F_\veps(\varphi) :=
\ioe\!\! f\varphi\;dx +
b\left(\frac{\veps}{r_\veps}\right)^3\!\!\!
\ise  g\varphi dx.
\ee
Then, for $\a>0$ (we shall choose a suitable value for this parameter later),
we can present the variational  formulation of
the problem (\ref{pbe1})--(\ref{pbe8}):

\hspace{1.5cm}
To find
$(u^\veps,\theta^\veps)\in V_\veps\times H^1_0(\Om)$ such that
\be
\lab{guv}
  \forall (v,\varphi)\in V_\veps\times H^1_0(\Om),\quad \la
G(u^\veps,\theta^\veps),(v,\varphi)\ra =
\F_\veps(\varphi)
\ee
where the mapping $G:\,V_\veps\times H^1_0(\Om)\to V_\veps'\times
H^{-1}(\Om)$ is defined by
$$
\ba{lll}
  \la G(u,\theta),(v,\varphi)\ra \!\!\!\!\!&= &\!\!\!
\dis\a\ioe \!\!\na u\na v\;dx - \a a \ioe \theta v_3\;dx
\\
&&+
\!\!\dis\ioe \!\!\!\na\theta\na\varphi\;dx + \ioe u \varphi\na\theta
\;dx + b\!\!\left(\frac{\veps}{r_\veps}\right)^3
\dis\!\!\!\!\int_{T_\veps}\!\!\!\na\theta\na\varphi\;dx.
\ea
$$

In order to prove the existence theorem for problem
(\ref{guv}), we make use of the following result of
Gossez.

\bthm
Let  $X$ be a reflexive Banach space and $G:\, X\to X'$ a continuous
mapping between the corresponding weak topologies. If
$$
  \frac{\la G\varphi,\varphi\ra}{\vert\varphi\vert_X}\to
\infty\quad\mbox{as}\quad \vert\varphi\vert_X\to\infty
$$
then $G$ is a surjection.
\ethm

Acting as in the proof of Theorem~5.2.2~\cite{EneP} Ch~1, Sec.~5, we find
that the existence of the weak solutions of problem (\ref{guv}) is
assured if  $\a$ is chosen  sufficiently small.

Moreover, if
$(u^\veps,\theta^\veps)$ is a solution of problem
(\ref{guv}), then, by using the weak maximum principle, we
obtain that
$\theta^\veps\in L^\infty(\Om)$, (see Theorem~3.4~\cite{EneP} Ch~2,
Sec.~3).

\brem
{\rm
For any $a>0$, we have proved the existence of a solution of
(\ref{guv}), but we do not have a uniqueness result,
except if we assume that $a>0$ is small enough.
}
\erem

In the sequel, $C$ will denote a suitable positive constant independent
of $\veps$ and which may differ from line to line.

\section{Basic inequalities}\lab{basic}

Lemma~\ref{l:dr1r2} and Lemma~\ref{l:ar} below are set without proof
since it is an adaptation of the case $p=2$ of Lemma~A.3~\cite{BellieudB}
and Lemma~A.4~\cite{BellieudB} respectively but with integrals set on
spheres.

\blem
\lab{l:dr1r2}
For every   $0 <r_1<r_2$, consider:
$$
  C(r_1,r_2) := \{x\in\R^3,\quad r_1<\vert x\vert < r_2\}.
$$
Then, if $u\in H^1(C(r_1,r_2))$, the following estimate holds true:
\be
\lab{c:dr1r2}
 \vert\na u\vert^2_{C(r_1,r_2)}
\geq  \frac{4\pi r_1r_2}{r_2 - r_1}\left\vert
\intb_{\Sp_{r_2}}u\;d\sigma -
\intb_{\Sp_{r_1}}u\;d\sigma\right\vert^2,
\ee
where
$$
  \intb_{\Sp_r}\cdot\;d\sigma :=
\frac{1}{4\pi r^2}\int_{\Sp_r}\cdot\;d\sigma.
$$

\elem

\blem
\lab{l:ar}
There exists a positive constant $C>0$ such that: $\forall (R,\a)\in
\R^+\times (0,1)$, $\forall u\in H^1(B(0,R))$,
$$
  \int_{B(0,R)}\vert u- \intb_{\Sp_{\a R}}u\;d\sigma\vert^2\;dx\leq
C\frac{R^2}{\a}\vert\na u\vert^2_{B(0,R)}.
$$
\elem

From now on, we denote by $R_\veps$ a radius with the property
$r_\veps<< R_\veps <<\veps$, that is :

\be \lab{re} \quad \lim_{\veps\to 0}\frac{r_\veps}{R_\veps} =
\lim_{\veps\to 0} \frac{R_\veps}{\veps} = 0. \ee

\noindent Obviously, its existence is insured by the assumption $0
< r_\veps <<\veps$.

We introduce the measure
$$
  dm^\veps :=
\frac{3}{4\pi}\left(\frac{\veps}{r_\veps}\right)^31_{T_\veps}(x)\;dx
$$
and denote the norm in $L^2_{\me}$ by:
$$
  \vert\varphi\vert^2_\me := \int\vert\varphi\vert^2\;d\me.
$$

We denote the domain confined between the spheres of
radius
$a$ and
$b$ by
$$
  \C(a,b) :=
\{x\in\R^3,\; a < \vert x\vert< b\}
$$
and correspondingly
$$
  \C^k(a,b) := \veps  k + \C(a,b),
$$

We also use the following notations:
$$
  \C_\veps := \cup_{k\in\Z_\veps} \C^k(r_\veps,R_\veps).
$$
$$
  S_{r_\veps}^k = \pa B(\veps k,r_\veps),\quad
S_{r_\veps} := \cup_{k\in\Z_\veps}S^k_{r_\veps},
$$
$$
  S_{R_\veps}^k = \pa B(\veps k,R_\veps),\quad
S_{R_\veps} := \cup_{k\in\Z_\veps}S^k_{R_\veps},
$$

Consider the piecewise constant functions defined
after some $\theta \in H^1_0(\Om)$ by
\beq
\lab{ttau}
\ti\tau^\veps(x)
& = &
\sum_{k\in\Z_\veps}\left(\intb_{\Sp^k_{r_\veps}}\theta\;d\sigma
\right) 1_{Y^k_\veps}(x),
\\
\lab{ttheta}
\ti\theta^\veps(x)
& = &
\sum_{k\in\Z_\veps}\left(\intb_{\Sp^k_{R_\veps}}\theta\;d\sigma
\right) 1_{Y^k_\veps}(x).
\eeq

\blem
For every $\theta\in H^1_0(\Om)$, we have
\beq
\lab{I1}
  \io\vert\theta - \ti\theta^\veps\vert^2\;dx
& \leq &
C\frac{\veps^3}{R_\veps}\io\vert\na\theta\vert^2\;dx,
\\
\lab{I2}
\ise\vert\theta - \ti\tau^\veps\vert^2\;dx
& \leq &
C r_\veps^2\ise\vert\na\theta\vert^2\;dx
\\
\lab{I3}
\io\vert\ti\theta^\veps - \ti\tau^\veps\vert^2\;dx
&\leq &
C\frac{\veps^3}{r_\veps}\ibe\vert\na\theta\vert^2\;dx.
\eeq
where $\ti\theta^\veps$ and $\ti\tau^\veps$ are defined by (\ref{ttau})
and (\ref{ttheta}).

Moreover:
\be
\lab{mexb}
  \io\vert\ti\theta^\veps\vert^2\;dx =
\int\vert\ti\theta^\veps\vert^2\;d\me,\quad
  \io\vert\ti\tau^\veps\vert^2\;dx =
\int\vert\ti\tau^\veps\vert^2\;d\me.
\ee

\elem
\bP
Notice that by definition:
$$
\io\vert\theta - \ti\theta^\veps\vert^2\;dx
=
\sum_{k\in\Z_\veps}\iyke\vert \theta -
\intb_{\Sp^k_{R_\veps}}
\theta\;d\sigma\vert^2\;dx
\leq\sum_{k\in\Z_\veps}\int_{B(\veps k,\frac{\veps\sqrt{3}}{2})}\vert
\theta -
\intb_{\Sp^k_{R_\veps}}
\theta\;d\sigma\vert^2\;dx
$$
where we have used that
$$
  Y^k_\veps\subset B(\veps k,\frac{\veps\sqrt{3}}{2})
$$
for every $k\in\Z_\veps$. We use Lemma~\ref{l:ar} with
$$
  R = \frac{\veps\sqrt{3}}{2},\quad \a = \frac{2R_\veps}{\veps\sqrt{3}}
$$
to deduce that
$$
\io\vert\theta - \ti\theta^\veps\vert^2\;dx
\leq
C\left(\frac{\veps\sqrt{3}}{2}\right)^2\frac{\veps\sqrt{3}}{2R_\veps}
\sum_{k\in\Z_\veps}
\int_{B(\veps k,\frac{\veps\sqrt{3}}{2})}\vert\na\theta\vert^2\;dx
$$
$$
\leq
C\frac{\veps^3}{R_\veps}\sum_{k\in\Z_\veps}
\int_{B(\veps k,\frac{\veps\sqrt{3}}{2})}\vert\na\theta\vert^2\;dx
\leq
C\frac{\veps^3}{R_\veps}\io\vert\na\theta\vert^2\;dx
$$
which shows (\ref{I1}).

To establish (\ref{I2}), we recall the definition:
$$
\ise\vert\theta - \ti\tau^\veps\vert^2\;dx
=
\sum_{k\in\Z_\veps}\int_{B(\veps k,r_\veps)}\vert \theta -
\intb_{\Sp^k_{r_\veps}}\theta\;d\sigma\vert^2\;dx
$$
Applying Lemma~\ref{l:ar} with $R = r_\veps$ and $\a = 1$,
we get the result
\beqs
\ise\vert\theta - \ti\tau^\veps\vert^2\;dx
& \leq &
Cr_\veps^2\sum_{k\in\Z_\veps}\int_{B(\veps k,r_\veps)}
\vert\na\theta\vert^2\;dx
\leq C r_\veps^2\ise\vert\na\theta\vert^2\;dx.
\eeqs

We come to (\ref{I3}). Indeed, applying Lemma~\ref{l:dr1r2} and
(\ref{re}):
$$
\io\vert\ti\theta^\veps - \ti\tau^\veps\vert^2\;dx
=
\sum_{k\in\Z_\veps}\iyke\vert\intb_{\Sp^k_{R_\veps}}
\theta\;d\sigma
- \intb_{\Sp^k_{r_\veps}}\theta\;d\sigma\vert^2\;dy
$$
$$
\leq
\sum_{k\in\Z_\veps}\iyke\frac{(R_\veps - r_\veps)}{4\pi R_\veps r_\veps}
\;dy\int_{C^k_{r_\veps,R_\veps}}
\vert\na\theta\vert^2\;dx
=
\frac{(R_\veps - r_\veps)}{4\pi r_\veps R_\veps}
\sum_{k\in\Z_\veps}\veps^3
\int_{C^k_{r_\veps,R_\veps}}\vert\na\theta\vert^2\;dx
$$
$$
=
C\veps^3\frac{(R_\veps - r_\veps)}{4\pi r_\veps R_\veps}\ibe
\vert\na\theta\vert^2\;dx
\leq
C\frac{\veps^3}{r_\veps}\ibe\vert\na\theta\vert^2\;dx.
$$
Finally, a direct computation yields (\ref{mexb}).
\eP

\bprop
\lab{p:naom}
For any $\theta\in H^1_0(\Om)$, there holds true:
$$
    \int\vert\theta\vert^2\;d\me\leq
C\max{(1,\frac{\veps^3}{r_\veps})}
\io\vert\na\theta\vert^2\;dx.
$$
\eprop
\bP
We have:
$$
\int\vert\theta\vert^2\;d\me
\leq
2\int\vert\theta - \ti\tau^\veps\vert^2\;d\me +
2 \int\vert\ti\tau^\veps\vert^2\;d\me
$$
$$
=
2\int\vert\theta - \ti\tau^\veps\vert^2\;d\me +
2 \int_\Om\vert\ti\tau^\veps\vert^2\;dx
$$
$$
\leq
Cr^2_\veps\int\vert\na\theta\vert^2\;d\me +
4 \int_\Om\vert\ti\tau^\veps - \ti\theta^\veps\vert^2\;dx +
8 \int_\Om\vert\ti\theta^\veps - \theta\vert^2\;dx +
8 \int_\Om\vert\theta\vert^2\;dx
$$
$$
\leq
Cr^2_\veps\left(\frac{\veps}{r_\veps}\right)^3
\ise\vert\na\theta\vert^2\;dx
+  C \frac{\veps^3}{r_\veps}
\int_{\C_\veps}\vert\na\theta\vert^2\;dx +
C\frac{\veps^3}{R_\veps}\int_\Om\vert\na\theta\vert^2\;dx
+
C\io\vert\na\theta\vert^2\;dx
$$
$$
\leq
C\left(\frac{\veps^3}{r_\veps} + \frac{\veps^3}{R_\veps} + 1\right)
\io\vert\na\theta\vert^2\;dx
\leq C\max{(1,\frac{\veps^3}{r_\veps})}\io\vert\na\theta\vert^2\;dx
$$
\eP

\blem
\lab{l:ffe}
For $\varphi\in C_c(\Om)$ consider the piecewise constant function:
$$
  \varphi^\veps(x) :=
\sum_{k\in\Z_\veps}\left(\intb_{Y^k_\veps}\varphi\;dx
\right)\;1_{B(\veps k,r_\veps)}(x).
$$
Then:
$$
  \lim_{\veps\to 0}\vert\varphi - \varphi^\veps\vert_\me = 0.
$$
\elem
\bP
Notice that
$$
  \vert\varphi - \varphi^\veps\vert^2_\me =
\frac{3}{4\pi}\left(\frac{\veps}{r_\veps}\right)^3\sum_{k\in\Z_\veps}
\int_{B(\veps k,r_\veps)}\vert\varphi - \intb_{Y^k_\veps}\varphi\;dy
\vert^2\;dx.
$$
As we have also
$$
  \vert B(\veps k,r_\veps)\vert = \frac{4\pi}{3}r^3_\veps,\quad
{\rm card}(\Z_\veps) \simeq \frac{\vert\Om\vert}{\veps^3}
$$
then, by the uniform continuity of $\varphi$ on $\Om$, the result follows.
\eP

\section{A priori estimates}

In the sequel, we denote
\be
\lab{gam}
  \gamma_\veps := \frac{r_\veps}{\veps^3}
\ee
and  we assume that
\be
\lab{cap}
 \lim_{\veps\to 0}\gamma_\veps =  \gamma\in ]0,+\infty[.
\ee

We denote $\F\in H^{-1}(\Om)$ by
\be
\lab{F}
 \F(\varphi) := \io
f\varphi\;dx +
\frac{4\pi b}{3}\io g\varphi\;dx
\ee

\bprop
\lab{p:3.1}
We have
$$
\F_\veps\wto\F\quad\mbox{weakly in}\quad H^{-1}(\Om)
$$
\eprop
\bP
For $\varphi\in H^1_0(\Om)$ it follows
$$
  \vert\F_\veps(\varphi)\vert\leq
\vert f\vert_\Ome\vert\varphi\vert_\Ome +
C\left(\frac{\veps}{r_\veps}\right)^3
\vert g\vert_\infty\left\vert
\ise \varphi\;dx\right\vert
$$
\be
\lab{estf}
  \leq C\vert\varphi\vert_\Om + C\left\vert\int \varphi\;d\me\right\vert
\ee
with
\be
\lab{estf1}
  \left\vert\int \varphi\;d\me\right\vert
\leq
(\int\;d\me)^{1/2}
\left(\int\vert\theta^\veps\vert^2\;d\me\right)^{1/2}=
\sqrt{\vert\Om\vert}\left(\int\vert\varphi\vert^2\;d\me\right)^{1/2}.
\ee
Notice that due to (\ref{cap}), Proposition~\ref{p:naom} also reads
\be
\lab{tme1}
  \int\vert\varphi\vert^2\;d\me\leq C\vert\na\varphi\vert^2_\Om.
\ee
Substituting (\ref{estf1}) and (\ref{tme1}) into the right-hand side of
(\ref{estf}), we get, using Poincar\'e's inequality,
\be
\lab{estf3}
  \vert\F_\veps(\varphi)\vert\leq C\vert\na\varphi\vert_\Om.
\ee
Now, let $\varphi\in \DD(\Om)$.
By the Mean Theorem, there exist
$\xi^k_\veps\in B(\veps k,r_\veps)$ such that
\beqs
\F_\veps(\varphi)
& = &
\ioe f\varphi\;dx +
b\left(\frac{\veps}{r_\veps}\right)^3\sum_{k\in\Z^\veps}\int_{B(\veps
k,r_\veps)} g(x)\varphi(x)\;dx
\\
& = &
\ioe f\varphi\;dx +
b\left(\frac{\veps}{r_\veps}\right)^3\sum_{k\in\Z^\veps}
\frac{4\pi}{3}r^3_\veps g(\xi^k_\veps)\varphi(\xi^k_\veps)
\\
& = &
\ioe f\varphi\;dx +
\frac{4\pi b}{3}
\sum_{k\in\Z^\veps}\vert Y^k_\veps\vert
g(\xi^k_\veps)\varphi(\xi^k_\veps).
\eeqs

There follows
\be
\lab{F}
  \forall\varphi\in \DD(\Om),\quad\lim_{\veps\to 0}\F_\veps(\varphi) = \io
f\varphi\;dx +
\frac{4\pi b}{3}\io g\varphi\;dx = \F(\varphi).
\ee
The proof is completed by (\ref{estf3}) and the density of $\DD(\Om)$ in
$H^1_0(\Om)$.
\eP

\bprop
If $( u^\veps,\theta^\veps)\in V_\veps\times H^1_0(\Om)$ is a solution of
the problem (\ref{guv}), and if $\hat u^\veps$ stands
for $u^\veps$ continued with zero to $\Om$, then we have
\be
\lab{ut}
  \hat u^\veps\quad\mbox{and}\quad
\theta^\veps\quad\mbox{are bounded in $H^1_0(\Om)$}.
\ee
Moreover,
\be
\lab{c:qme}
  \vert\na\theta^\veps\vert^2_\Ome  +
b\left(\frac{\veps}{r_\veps}\right)^3
\vert\na\theta^\veps\vert^2_{T_\veps}\leq C.
\ee
\eprop
\bP
Substituting $v = u^\veps$ in
(\ref{vare1}) and noticing that
$$
    \ioe u^\veps\theta^\veps \na\theta^\veps\!dx=
\ioe\!\!\!u^\veps\na\!\!\left(\frac{\vert\theta^\veps\vert^2}{2}\right)\!dx
=
-\ioe\!\!\!\div (u^\veps)\!\!
\left(\frac{\vert\theta^\veps\vert^2}{2}\right)\!dx = 0,
$$
we get:
\beq
\lab{est1}
\vert\na u^\veps\vert_\Ome & \leq &
a\vert\theta^\veps\vert_\Ome,
\eeq
Seting $\varphi = \theta^\veps$ in (\ref{vare4}) and taking into account
Proposition~\ref{p:3.1}, we find
\be
\lab{est2}
  \vert\na\theta^\veps\vert^2_\Ome  +
b\left(\frac{\veps}{r_\veps}\right)^3
\vert\na\theta^\veps\vert^2_{T_\veps} = \F_\veps(\theta^\veps)\leq
C\vert\na\theta^\veps\vert_\Om
\ee

Noticing that
$b\left(\frac{\veps}{r_\veps}\right)^3>>1$, we deduce from (\ref{est2}):
$$
  \vert\na\theta^\veps\vert_\Om^2\leq
\vert\na\theta^\veps\vert^2_\Ome  +
b\left(\frac{\veps}{r_\veps}\right)^3
\vert\na\theta^\veps\vert^2_{T_\veps}
\leq
C\vert\na\theta^\veps\vert_\Om.
$$
Therefore
\be
\lab{c:name}
  \vert\na\theta^\veps\vert_\Om\leq C
\ee
and thus
\be
\lab{ta}
  \vert\theta^\veps\vert_\Om\leq C.
\ee
Then, (\ref{c:qme}) follows from (\ref{est2}).
Finally, (\ref{ut}) is completed by the estimates (\ref{est1})
and (\ref{ta}).
\eP

\bprop
\lab{p:uttau}
There exist $ u\in H^1_0(\Om;\R^3)$, $\theta\in H^1_0(\Om)$ and
$\tau\in L^2(\Om)$ such that, on some subsequence,
\beqs
 \hat u^\veps & \wto & u\quad\mbox{in}\quad H^1_0(\Om;\R^3),
\\
\theta^\veps & \wto & \theta\quad\mbox{in}\quad H^1_0(\Om),
\\
\ti\tau^\veps & \wto & \tau\quad\mbox{in}\quad L^2(\Om),
\\
\theta^\veps\;d\me & \buildrel{\star}\over{\wto} &
\tau\;dx\quad\mbox{in}\quad
\M_b(\ov\Om),
\eeqs
where $\M_b(\ov\Om)$ is the set of bounded Radon measures on $\ov\Om$ and
where $\buildrel{\star}\over{\wto}$ denotes the weak-star convergence in the
measures.
\eprop
\bP
From (\ref{ut}), we get, on some subsequence, the following convergences:
\be
\lab{teh}
  \theta^\veps  \wto \theta\quad\mbox{in}\quad H^1_0(\Om)
\ee
\be
\lab{qf}
  \theta^\veps\to\theta\quad\mbox{in}\quad L^2(\Om).
\ee
\be
\lab{ueu}
  \hat u^\veps\wto u\quad\mbox{in}\quad H^1_0(\Om;\R^3).
\ee
Moreover, (\ref{I1}) yields
$$
  \vert\theta^\veps - \ti\theta^\veps\vert^2_\Om\leq
C\frac{r_\veps}{\veps^3}\frac{r_\veps}{R_\veps}
\vert\na\theta^\veps\vert^2_\Om
$$
which obviously yields
$$
  \lim_{\veps\to 0}\vert\theta^\veps - \ti\theta^\veps\vert^2_\Om = 0.
$$
Combining with (\ref{qf}), we infer that
\be
\lab{ttild}
  \ti\theta^\veps\to\theta\quad\mbox{in}\quad L^2(\Om).
\ee
We set
\be
\lab{te}
  \tau^\veps :=
\frac{3}{4\pi}\left(\frac{\veps}{r_\veps}\right)^3\theta^\veps
1_{T_\veps}(x),
\ee
and hence
$$
  \theta^\veps\;d\me = \tau^\veps\;dx.
$$
Taking (\ref{ut}) and (\ref{tme1}) into account, we obtain

$$
  \int \vert\theta^\veps\vert^2\;d\me\leq  C.
$$
We also remark that for any $\varphi\in C_c(\Om)$, we have
$$
  \int\varphi\;d\me\to \io\varphi dx.
$$
Then, using Lemma A-2 of \cite{BellieudB}, we find that there
exists some $\tau\in L^2(\Om)$ such that, on some subsequence, the following
convergence holds:
\be
\lab{qme}
  \theta^\veps\;d\me\buildrel{\star}\over{\wto}
\tau\;dx,\quad
\M_b(\ov\Om).
\ee

Moreover, recall that from (\ref{I2}) we have, taking into account
(\ref{c:qme}):
\be
\lab{qtt}
  \int\vert\theta^\veps - \ti\tau^\veps\vert^2\;d\me\leq C r^2_\veps
\int\vert\na\theta^\veps\vert^2\;d\me\leq C r^2_\veps.
\ee
This implies:
$$
  (\theta^\veps - \ti\tau^\veps)\;d\me
\buildrel{\star}\over{\wto} h\;dx,\quad \M_b(\ov\Om)
$$
for some $h\in L^2(\Om)$ and
$$
  \vert h\vert^2_\Om\leq \liminf_{\veps\to 0}
\int\vert\theta^\veps - \ti\tau^\veps\vert^2\;d\me = 0,
$$
that is:
\be
\lab{tte}
  (\theta^\veps - \ti\tau^\veps)\;d\me
\buildrel{\star}\over{\wto} 0,\quad \M_b(\ov\Om).
\ee
Notice that from (\ref{I3}):
\be
\lab{ttc}
  \vert\ti\tau^\veps\vert^2_\Om
\leq
2\vert\ti\tau^\veps - \ti\theta^\veps\vert^2_\Om +
2\vert\ti\theta^\veps\vert^2_\Om
\leq
C\frac{\veps^3}{r_\veps}\vert\na\theta^\veps\vert^2_{\C_\veps} + C
\leq C,
\ee
and hence, for some $\ti\tau\in  L^2(\Om)$,
\be
\lab{tt}
  \ti\tau^\veps\wto\ti\tau\quad\mbox{in}\quad L^2(\Om).
\ee
Combining (\ref{qme}) and (\ref{tte}), we
arrive at
$$
  \ti\tau^\veps\;d\me\buildrel{\star}\over{\wto} \tau\;dx,\quad
\M_b(\ov\Om).
$$
It remains to show that
\be
\lab{tg}
  \ti\tau = \tau.
\ee
To that aim, let $\varphi\in C_c(\Om)$ and let
$$
  \varphi^\veps(x) :=
\sum_{k\in\Z_\veps}\left(\intb_{Y^k_\veps}\varphi\;dy\right)
\;1_{B(\veps k,r_\veps)}(x).
$$
We have
\beqs
&&
\vert\io (\tau^\veps - \ti\tau^\veps)\,\varphi\;dx\vert
=
\left\vert
\frac{3}{4\pi}
\left(\frac{\veps}{r_\veps}\right)^3\!\!\!\ise\theta^\veps\varphi\;dx
\!-\!\!\!
\io
\sum_{k\in\Z_\veps}\!\!\left(\intb_{\Sp^k_{r_\veps}}\theta^\veps\;d\sigma
\right)\!\!1_{Y^k_{\veps}}\varphi dx\right\vert
\\
&=&
\left\vert
\int\theta^\veps\varphi\;d\me -
\veps^3\sum_{k\in\Z_\veps}\left(\intb_{\Sp^k_{r_\veps}}\theta^\veps\;d\sigma
\right)
\intb_{Y^k_\veps}\varphi\;dx\right\vert
=
\left\vert
\int\theta^\veps\varphi\;d\me -
\int\ti\tau^\veps\varphi^\veps\;d\me
\right\vert
\\
&& \hspace{2.5cm} \leq
\left\vert\io (\theta^\veps - \ti\tau^\veps)\varphi\;d\me \right\vert +
\left\vert \io\ti\tau^\veps(\varphi - \varphi^\veps)\;d\me\right\vert
\eeqs
\beq
\lab{est3}
& \leq &
\vert\theta^\veps -\ti\tau^\veps\vert_\me \vert\varphi\vert_\me +
\vert\ti\tau^\veps\vert_\me\vert\varphi - \varphi^\veps\vert_\me.
\eeq
From (\ref{mexb}) and (\ref{ttc}), we deduce that
$$
  \vert\ti\tau^\veps\vert_\me = \vert\ti\tau^\veps\vert_\Om\leq C.
$$
Moreover, $\varphi\in C_c(\Om)$ yields
$$
  \vert\varphi\vert_\me\leq C.
$$
Then, (\ref{est3}) becomes
$$
  \vert\io (\tau^\veps - \ti\tau^\veps)\,\varphi\;dx\vert\leq
C\vert\theta^\veps -\ti\tau^\veps\vert_\me  +
C\vert\varphi - \varphi^\veps\vert_\me.
$$
From (\ref{qtt}), we infer that
\be
\lab{est4}
  \vert\io (\tau^\veps - \ti\tau^\veps)\,\varphi\;dx\vert\leq
Cr_\veps +
C\vert\varphi - \varphi^\veps\vert_\me.
\ee
Thus (\ref{est4}) and Lemma~\ref{l:ffe} yield
$$
  \lim_{\veps\to 0}\io (\tau^\veps - \ti\tau^\veps)\,\varphi\;dx = 0.
$$
As this holds for every $\varphi\in C_c(\Om)$, the density of
$C_c(\Om)$ in $L^2(\Om)$ together with (\ref{tt}) and (\ref{qme})
imply that
$
  \tau^\veps\wto\ti\tau = \tau\quad\mbox{in}\quad L^2(\Om).
$
\eP

\section{The two macroscopic heat equations}

The aim of this section is to pass to the limit as $\veps\to 0$ in the
variational formulation
\beq
\lab{vare}
\ba{lll}
  \forall\Phi\in H^1_0(\Om),
\!\!\!\!&&\!\!\!\!
\dis\ioe\na\theta^\veps\na\Phi\;dx +
b\left(\frac{\veps}{r_\veps}\right)^3
\ise \na\theta^\veps
\na\Phi\;dx +
\\
&& + \dis\ioe u^\veps\na\theta^\veps\Phi\;dx
=
\F_\veps(\Phi).
\ea
\eeq

Let $\varphi,\psi\in \DD(\Om)$ and set
\beq
\lab{fie}
\varphi^\veps(x)
& = &
\sum_{k\in\Z_\veps}\left(\intb_{\Sp^k_{r_\veps}}\varphi\;d\sigma\right)
1_{Y^k_\veps}(x),
\\
\lab{psie}
\psi^\veps(x)
& = &
\sum_{k\in\Z_\veps}\left(\intb_{\Sp^k_{r_\veps}}\psi\;d\sigma\right)
1_{Y^k_\veps}(x).
\eeq
Let $W^\veps$ denote the fundamental solution of the
Laplacian, namely
\beq
\lab{defbwe1}
  \Delta W^\veps & = & 0\quad\mbox{in}\quad \C(r_\veps,R_\veps),
\\
W^\veps & = & 1\quad\mbox{in}\quad  r = r_\veps,
\\
\lab{defbwe2}
W^\veps & = & 0\quad\mbox{in}\quad  r = R_\veps.
\eeq
The same arguments as in the proof of Lemma~A.3~\cite{BellieudB} yield
\be
\lab{defbwe}
  W^\veps(r) =  \frac{r_\veps}{(R_\veps -
r_\veps)}\left(\frac{R_\veps}{r} - 1\right)\quad\mbox{if}\quad
y\in \C(r_\veps,R_\veps)\quad\mbox{and}\quad \vert y\vert = r.
\ee
Then, we set
\beq
\lab{defwe}
  w^\veps(x) &:=& \left\{\ba{c}
0\quad\mbox{in}\quad\Ome\sm \C_\veps,
\\
W^\veps( x - \veps k)\quad\mbox{in}\quad  \C^k_\veps,\quad
\forall k\in\Z_\veps,
\\
1\quad\mbox{in}\quad T_\veps.
\ea\right.
\eeq

\bprop
\lab{p:web}
We have
\be
\lab{web}
  \vert\na w^\veps\vert_\Om\leq C
\ee
\eprop
\bP
Indeed, direct computation shows
\beqs
\vert\na w^\veps\vert_\Om^2
& = &
\sum_{k\in\Z_\veps}\int_{C^k_{r_\veps,R_\veps}}
\vert\na w^\veps\vert^2\;dx
\\
& = &
\sum_{k\in\Z_\veps}\int_0^{2\pi}\;d\Phi\int_0^{\pi}\sin\Theta\;d\Theta
\int_{r_\veps}^{R_\veps}\frac{dr}{r^2}\left(\frac{r_\veps
R_\veps}{R_\veps - r_\veps}\right)^2
\\
& \leq &
C\frac{\vert\Om\vert}{\veps^3}
\left(\frac{1}{r_\veps} - \frac{1}{R_\veps}\right)\left(\frac{r_\veps
R_\veps}{R_\veps - r_\veps}\right)^2
\leq C\frac{\gamma_\veps}{(1 - \frac{r_\veps}{R_\veps})}.
\eeqs
The proof is completed by (\ref{re}) and (\ref{cap}).
\eP

For $\varphi,\psi\in\DD(\Om)$, let us define
\be
\lab{deffie}
\Phi^\veps = (1-w^\veps)\varphi + w^\veps\psi^\veps.
\ee

\blem
\lab{p:fie}
We have
$$
  \lim_{\veps\to 0}\vert\Phi^\veps - \varphi\vert_\Om = 0.
$$
\elem
\bP
First notice that $w^\veps\to 0$ in $L^2(\Om)$. Indeed:
$$
  \vert w^\veps\vert_\Om = \vert w^\veps\vert_{\C_\veps\cup
T_\veps}\leq \vert \C_\veps\cup T_\veps\vert =
\frac{\vert\Om\vert}{\veps^3}
\frac{4\pi}{3}R_\veps^3
$$
and $\lim_{\veps\to 0}\frac{R_\veps}{\veps} = 0$ by assumption (\ref{re}).
As an immediate consequence:
$$
  (1-w^\veps)\varphi \to \varphi\quad\mbox{in}\quad L^2(\Om).
$$
Moreover, the uniform continuity of $\psi$ over $\Om$
implies that
$$
  \lim_{\veps\to 0}\vert\psi^\veps - \psi\vert_\infty = 0
$$
so that
$$
  w^\veps\psi^\veps = w^\veps(\psi^\veps - \psi) + w^\veps\psi\to 0
\quad\mbox{in}\quad L^2(\Om).
$$
This achieves the proof.
\eP

\bprop
\lab{p:vare}
If $\theta^\veps$ is solution of (\ref{guv}) and $\Phi^\veps$ is given
by (\ref{deffie}) for any $\psi,\varphi\in\DD(\Om)$, then  we have
\beqs
&&
  \lim_{\veps\to 0}\ioe \na\theta^\veps\cdot\left(
\na\Phi^\veps + \Phi^\veps u^\veps\right)\;dx
\\
& =&
\io \na\theta\cdot\left(\na\varphi + \varphi u\right)\;dx +
  4\pi\gamma
\io(\theta - \tau)(\psi -\varphi)\;dx.
\eeqs
\eprop
\bP
First consider%
$$
  \ioesmbe \na\theta^\veps\cdot\left(
\na\Phi^\veps + \Phi^\veps u^\veps\right)\;dx
$$
which reduces to
$$
  \ioesmbe \na\theta^\veps\cdot\left(
\na\varphi + \varphi u^\veps\right)\;dx =
\io  \na\theta^\veps\cdot\left(
\na\varphi 1_{\Ome\sm \C_\veps} + \varphi 1_{\Ome\sm
\C_\veps} u^\veps\right)\;dx.
$$
Lebesgue's dominated convergence theorem yields
$
  \na\varphi 1_{\Ome\sm \C_\veps}\to \na\varphi
$
in $L^2(\Om)$. Thus, taking (\ref{teh}) into account:
$$
  \io\na\theta^\veps\cdot\na\varphi 1_{\Ome\sm
\C_\veps}\;dx\to
\io\na\theta\cdot\na\varphi\;dx.
$$
Moreover,
$$
  \vert 1_{\Ome\sm \C_\veps} u^\veps -  u\vert_\Om\leq \vert
u^\veps -  u\vert_\Om +
\vert u\vert_{\C_\veps\cup T_\veps}
$$
and the right-hand side converges to zero because
(\ref{ueu}) yields
\be
\lab{uef}
  u^\veps\to  u\quad\mbox{in}\quad L^2(\Om)
\ee
and we apply  Lebesgue's
dominated convergence theorem to conclude with the second term. Thus
\be
\lab{ue}
   1_{\Ome\sm \C_\veps} u^\veps \to  u\quad\mbox{in}\quad
L^2(\Om).
\ee
Now, as $\varphi\in C_c(\Om)$,
$
  \varphi 1_{\Ome\sm \C_\veps} u^\veps \to \varphi u
$
in $L^2(\Om)$. Thus, using (\ref{teh}) again,
$$
  \io  \na\theta^\veps\cdot \varphi 1_{\Ome\sm \C_\veps}
u^\veps\;dx\to
\io  \na\theta\cdot \varphi  u\;dx.
$$
As a result:
\be
\lab{mbe}
\lim_{\veps\to 0} \ioesmbe \na\theta^\veps\cdot\left(
\na\Phi^\veps + \Phi^\veps u^\veps\right)\;dx =
 \io\na\theta\cdot\left(
\na\varphi + \varphi u\right)\;dx.
\ee
Now, we come to the remaining part, namely
\be
\lab{be0}
\ba{lcl}
&&
  \dis\ibe \na\theta^\veps\cdot\left(
\dis\na\Phi^\veps + \Phi^\veps u^\veps\right)\;dx
=
  \ibe \na\theta^\veps\cdot\left(
\na\varphi + \varphi u^\veps\right)\;dx
\\
\!\!\!\!\!
&&\!\!\!\!\!+
\dis\ibe \!\!\na\theta^\veps\!\!\cdot\!\left(
\na w^\veps (\psi^\veps - \varphi) + w^\veps (-\na\varphi) +
 w^\veps u^\veps(\psi^\veps - \varphi)
\right) dx
\\
& := & I_1 + I_2
\ea
\ee
We have
\be
\lab{be1}
   I_1= \ibe \na\theta^\veps\cdot\left(
\na\varphi + \varphi u^\veps\right)\;dx =
 \ibe \na\theta^\veps\cdot\na\varphi\;dx +
 \ibe \na\theta^\veps\cdot\varphi u^\veps\;dx.
\ee
In the first term, $1_{\C_\veps}\na\varphi\to 0$ in $L^2(\Om)$
and
$\na\theta^\veps\wto\na\theta$ in $L^2(\Om)$ imply
\be
\lab{be2}
  \ibe \na\theta^\veps\cdot\na\varphi\;dx\to 0.
\ee
The second term in (\ref{be1}) is handled by using the estimate:
$$
 \vert u^\veps\vert_{\C_\veps} = \vert 1_{\C_\veps} u^\veps\vert_\Om
\leq
\vert  u^\veps -  u\vert_\Om +
\vert  u\vert_{\C_\veps},
$$
where the right hand side tends to zero due to (\ref{uef}). Using
$\na\theta^\veps\wto\na\theta$ in $L^2(\Om)$ again, we deduce that
\be
\lab{be3}
  \ibe \na\theta^\veps\cdot\varphi u^\veps\;dx\to 0,
\ee
and hence $I_1$ tends to zero.

It
remains to study the integral $I_2$ in (\ref{be0}).
To that aim, first notice that
\beqs
&&
  I_2 = \ibe \na\theta^\veps\cdot
\na w^\veps (\psi^\veps - \varphi)\;dx =
\eeqs
\beq
\lab{be5}
&=& \!\!\!
\ibe \!\!\!\!\na\theta^\veps\cdot
\na w^\veps (\psi^\veps - \varphi^\veps)\;dx
+\ibe \!\!\!\! \na\theta^\veps\cdot
\na w^\veps (\varphi^\veps - \varphi)\;dx
\eeq
where $\varphi^\veps$  has been defined by (\ref{fie}). The second term in the
right-hand side of (\ref{be5}) may be estimated by
\be
\lab{be6}
  \vert \ibe \na\theta^\veps\cdot
\na w^\veps (\varphi - \varphi^\veps)\;dx\vert\leq
\vert\na\theta^\veps\vert_\Om\vert\na w^\veps\vert_\Om\vert\varphi -
\varphi^\veps\vert_\infty.
\ee
As $(w_\veps)$ is bounded in $H^1(\Om)$,
({\it see} Proposition~\ref{p:web}),
the right hand side
of (\ref{be6}) tends to zero by the uniform continuity of
$\varphi$ over $\Om$.

Going back to the first term in
the right hand side of (\ref{be5}), we may write
\beqs
&&
\ibe \na\theta^\veps\cdot
\na w^\veps (\psi^\veps - \varphi^\veps)\;dx
\\
& = &
\sum_{k\in\Z_\veps}\int_0^{2\pi}d\Phi\!\!\int_{0}^{\pi}\sin\Theta\;d\Theta
\int_{r_\veps}^{R_\veps}\left.\frac{\pa\theta^\veps}{\pa
r}\right\vert_{\C^k(r_\veps,R_\veps)}\!
\frac{d\ov w_\veps}{dr} r^2\;dr
\left(\intb_{\Sp^k_{r_\veps}}\psi\;d\sigma-
\intb_{\Sp^k_{r_\veps}}\varphi\;d\sigma\right)
\\
& = &\!\!\!
\frac{r_\veps R_\veps}{(R_\veps - r_\veps)}
\sum_{k\in\Z_\veps}\int_{\Sp_1}(\theta^\veps\vert_{\vert x - \veps k\vert =
r_\veps} - \theta^\veps\vert_{\vert x - \veps k\vert = R_\veps})
\left(\intb_{\Sp^k_{r_\veps}}\psi\;d\sigma-
\intb_{\Sp^k_{r_\veps}}\varphi\;d\sigma\right)d\sigma_1
\\
& = &
\frac{4\pi r_\veps R_\veps}{\veps^3(R_\veps - r_\veps)}
\io (\ti\tau^\veps - \ti\theta^\veps)(\psi^\veps -\varphi^\veps)\;dx
= \frac{4\pi\gamma_\veps}{\left(1 - \frac{r_\veps}{R_\veps}\right)}
\io(\ti\tau^\veps - \ti\theta^\veps)(\psi^\veps -\varphi^\veps)\;dx
\eeqs
from which we infer that $I_2$ is converging to
$$
4\pi\gamma
\io(\tau - \theta)(\psi -\varphi)\;dx,
$$
and the proof is completed.
\eP

We are in the position to state a part of our  main result:

\bcor
\lab{c:lim}
The limit $(u,\theta,\tau)$ verifies the following equations:
\beq
\lab{pblim1}
u\na\theta -\Delta\theta + 4\pi\gamma (\theta - \tau) & = &
f\quad\mbox{in}\quad\Om,
\\
\lab{pblim2}
\gamma (\tau - \theta) & = &
\frac{ b}{3} g\quad\mbox{in}\quad\Om.
\eeq
\ecor
\bP
Consider the variational formulation (\ref{vare4}) with the
test function
$\Phi = \Phi^\veps$ defined by (\ref{deffie}) for any
$\varphi,\psi\in\DD(\Om)$. Then, the left-hand side tends to
\be
\lab{c:vare}
  \io \na\theta\cdot\left(\na\varphi + \varphi u\right)\;dx +
  4\pi\gamma
\io(\tau - \theta)(\psi -\varphi)\;dx.
\ee
This is a direct consequence of Proposition~\ref{p:vare} together with
the remark that
$$
  \int_{T_\veps}\na\theta^\veps\na\Phi^\veps\;dx = 0
$$
since $\Phi^\veps$ is constant on every $B(\veps k,r_\veps)$,
$k\in\Z_\veps$.

The convergence of the right-hand side is obtained by using the uniform
continuity of $\psi$ and by Proposition~\ref{p:3.1}. Thus we find
the variational formulation of
(\ref{pblim1})-(\ref{pblim2}) and the proof is completed.
\eP

\section{The homogenized problem}

Proposition~\ref{p:uttau} yields the existence of some $u\in
H^1_0(\Om;\R^3)$ with $\div{(u)} = 0$ and for which the following convergence
holds on some subsequence
$$
  \hat u^\veps\wto u\quad\mbox{in}\quad H^1_0(\Om;\R^3).
$$
From \cite{Allaire}, we find that there exists an extension of the pressure
(denoted by $\hat p^\veps$) and some $p\in L^2(\Om)$ such that
$$
  \hat p^\veps\wto p\quad\mbox{in}\quad L^2(\Om)/\R.
$$
We denote by $(w^k_\veps,q^k_\veps)\in H^1(\C(r_\veps,\frac{\veps}{2}))\times
L^2_0(\C(r_\veps,\frac{\veps}{2}))$ the only solution of the following Stokes
problem
\beqs
\div{w^k_\veps} & = & 0\quad\mbox{in}\quad\C(r_\veps,\frac{\veps}{2}),
\\
- \Delta w^k_\veps + \na q^k_\veps & = &
0\quad\mbox{in}\quad\C(r_\veps,\frac{\veps}{2}),
\\
w^k_\veps & = & 0\quad\mbox{if}\quad r = r_\veps,
\\
w^k_\veps & = & \e^{(k)}\quad\mbox{if}\quad r = \frac{\veps}{2}.
\eeqs
Consequently, we define
$$
  v^k_\veps(x) = \left\{\ba{ll}
0\quad\mbox{if}\quad x\in T_\veps,
\\
w^k_\veps\left( x - \veps i\right)
\quad\mbox{if}\quad x\in \C^i(r_\veps,\frac{\veps}2),\quad i\in\Z_\veps,
\\
\e^{(k)}\quad\mbox{if}\quad x\in
\Ome\sm\cup_{i\in\Z_\veps} \C^i(r_\veps,\frac{\veps}2).
\ea\right.
$$
For $\varphi\in\DD(\Om)$, we set $v = \varphi v^k_\veps$ in (\ref{vare1}) and
then using the energy method like in \cite{CioranescuM3} we find the equation
that the velocity field satisfies in $H^{-1}(\Om)$:
\be
\lab{M}
  -\Delta u + 6\pi\gamma u  = - \na p +
a\theta\e^{(3)}\quad\mbox{in}\quad
\Om.
\ee

Finally, we summarize the results of  Proposition~\ref{p:uttau},
Corollary~\ref{c:lim} together with the relation (\ref{M}) into our main
theorem.

\bthm If $(u^\veps, p^\veps)$ is a solution of problem
(\ref{guv}), then the following convergences hold on some
subsequence \beqs
 \hat u^\veps  &\wto&  u\quad\mbox{in}\quad H^1_0(\Om;\R^3),
\\
\theta^\veps & \wto & \theta\quad\mbox{in}\quad H^1_0(\Om),
\\
\theta^\veps\;d\me & \buildrel{\star}\over{\wto} &
\tau\;dx\quad\mbox{in}\quad \M_b(\ov\Om), \eeqs

\noindent where $(u,\theta)\in H^1_0(\Om;\R^3)\times H^1_0(\Om)$,
which stand for the macroscopic velocity and temperature of the
fluid, and $\tau\in L^2(\Om)$, which stands for the macroscopic
temperature of the vanished suspensions, form a solution of the
following system:
$$
  \div{u} = 0\quad\mbox{in}\quad \Om,
$$
$$
  -\Delta u + 6\pi\gamma u  = - \na p +
a\theta\e^{(3)}\quad\mbox{in}\quad \Om,
$$
$$
u\na\theta -\Delta\theta + 4\pi\gamma (\theta - \tau)=
f\quad\mbox{in}\quad\Om,
$$
$$
 4\pi\gamma (\tau - \theta)  =  \frac{ 4\pi b}{3}
g\quad\mbox{in}\quad\Om.
$$
\ethm

\brem {\rm  In the present case, with suspensions of critical
size, the Brinkman-Boussinesq equation was an expected result;
nevertheless, our proof is different from that of \cite{Allaire},
which treated the homogenization of the Navier-Stokes equations
for perforated domains in a similar case.} \erem

\brem {\rm Our two-temperature model, with $\gamma$ as transfer
coefficient, is the macroscopic effect of the assumption on the
the ratio of the fluid/solid conductivities.} \erem

\brem {\rm The appearance of the source term $\dis \frac{4\pi
b}{3} g$ in the second macroscopic heat equation is strictly the
consequence of the assumption on the microscopic radiation.}\erem

\vs5
{\bf Acknowledgements.} This work was done during the visit of Fadila Bentalha
and Dan Poli\c{s}evschi at the I.R.M.A.R.'s Department of Mechanics
(University of Rennes~1) whose support is gratefully acknowledged.


\newcommand{\noopsort}[1]{}

\vs5

*  University of Batna, Department of Mathematics, Batna, Algeria,

\vs5

 ** Universit\'e de Rennes1, I.R.M.A.R, Campus de Beaulieu,
35042 Rennes Cedex (France)

\vs5

*** I.M.A.R., P.O. Box 1-764, Bucharest (Romania).

\end{document}